\documentclass{mcom-l}

\usepackage{bm}
\usepackage{mathrsfs}
\usepackage{amssymb}
\usepackage{enumitem}
\usepackage{times}
\usepackage{pdfpages}
\usepackage{float}
\usepackage{subfig,graphicx}
\usepackage{csquotes}
\usepackage{hyperref}
\usepackage{amsmath}
\usepackage{ulem}
\usepackage{lipsum} 
\usepackage{nccmath}
\usepackage{bigints}
\usepackage{natbib}
\usepackage{multirow}
\usepackage{adjustbox}
\usepackage{threeparttable}
\usepackage{booktabs}
\usepackage{bbm}

\usepackage{pdfpages}
\usepackage{wrapfig}
\usepackage{color}

\usepackage{float}
\usepackage{tabularx}
\usepackage{caption}
\usepackage{booktabs}
\usepackage{cite}

\usepackage{graphicx}

\usepackage{algorithm}
\usepackage{algpseudocode}

\newtheorem{theorem}{Theorem}

\newtheorem{proposition}{Proposition}

\newcommand{\A}{{\mathcal A}}

\newcommand{\raw}{\rightarrow}

\newcommand{\nti}{n\rightarrow\infty}

\numberwithin{equation}{section}
\numberwithin{theorem}{section}
\numberwithin{corollary}{section}
\numberwithin{lemma}{section}
\numberwithin{proposition}{section}
\numberwithin{corollary}{section}
\numberwithin{example}{section}
\numberwithin{remark}{section}

\begin{document}

\title[Optimal dimension growth rates in hd CLT]{Optimal dimension growth rates in the high dimensional  Central Limit Theorems}

\author{Mayukh Choudhury}
\address{Department of Mathematics, Indian Institute of Technology Bombay, Mumbai 400076, India}
\email{214090002@iitb.ac.in}

\author{Debraj Das}
\address{Department of Mathematics, Indian Institute of Technology Bombay, Mumbai 400076, India}
\email{debrajdas@math.iitb.ac.in}





\begin{abstract}
We investigate the optimal dimension growth rates for the high-dimensional Central Limit Theorem (CLT) over half spaces, Euclidean balls and Borel convex sets, under the existence of fourth moments. 
It has been observed in the literature that $d=o(n^c)$ with $0< c\leq 1/2$ is sufficient for the CLT to hold over Euclidean balls. As a consequence, the CLT over half spaces should ideally hold under similar dimension growth rates. In case of Borel convex sets, however $d= o(n^{c})$ with $0< c\leq 2/5$ is generally sufficient. The most recent results for CLT over Euclidean balls and Borel convex sets with respective rates $d= o(\sqrt{n})$ and $d= o(n^{2/5})$ is due to \citet{fang2024large} and under the existence of fourth moments. 
In this paper, we investigate whether these sufficient growth rates of $d$ with respect to $n$ are also necessary. 
Thus our results completely characterize the sharpness of dimension growth rates under existence of fourth moments across half-spaces, Euclidean balls, and Borel convex sets.
\end{abstract}

\maketitle

\section{Introduction}\label{sec:intro4}

Central Limit Theorem (CLT) is one of the oldest results of probability theory which has vast applications in different fields. CLT is a universal statement on the distributional properties of sample mean. More precisely, it states that the distribution of a properly centered and scaled sample mean of a sequence of random vectors can be approximated by the standard Gaussian distribution. Although the 
form of the CLT was  known in the eighteenth century due to Abraham de Moivre, the modern theoretical developments around it only happened in the past hundred years, cf. \citet{fischer2011history}. The recent interest is in developing CLT uniformly over a collection of sets when the dimension $d$ of the underlying random vectors grows with the sample size $n$. A huge volume of work has already been done in this direction primarily motivated by the seminal works of \citet{senatov1981uniform}, \citet{gotze1991rate,talagrand1996transportation,bentkus2003dependence} and \citet{chernozhukov2013gaussian}; among others. In the fixed dimensional CLT, although the approximating Gaussian distribution depends only on the first two moments, the rate of convergence generally depends on the existence of moments higher than two. In high dimensional CLT, the rate of convergence translates to how fast $d$ can grow with $n$. It has been observed that the dimension $d$ can grow polynomially or exponentially with $n$ based on whether some polynomial or exponential moment of the underlying random vectors exist and also on the underlying collection of sets considered.  An important question in this regard is 
\begin{center}
``\textit{What is the optimal growth rate of the dimension $d$ as a function of the sample size $n$ for the validity of the CLT
in high dimensions?}''
\end{center}

Let $X_1,\dots, X_n$ be independent mean zero $d$-dimensional random vectors. Throughout this paper, we assume that the distribution of $X_i$ may depend on $n$ and $d$, unless specified otherwise. Define $T_n = \sum_{i=1}^{n}Y_i$ and assume that $\Sigma = Var(T_n)$ is non-singular. When $d$ is fixed then the Lindeberg-Feller CLT (cf. Theorem 11.1.1, \citet{athreya2006measure}) implies that the distribution of ${T}_n$ can be approximated by $Z$ where $Z$ has the Gaussian distribution with mean zero and covariance matrix $\Sigma$. The Berry-Esseen theorem (cf. Theorem 12.4 of \citet{bhattacharya1986normal}), on the other hand, tells us how good this approximation really is. However, the Berry-Esseen rate essentially depends on $d$ when it grows with $n$, and hence one can naturally ask what is the growth rate of $d$ above which the Gaussian approximation fails uniformly over a collection of sets. The specification of the collection of sets in the high dimensional CLT is important since the geometry of the collection has been observed to have significant effect on the Berry-Esseen rate. In particular, the literature broadly bifurcates into two geometric regimes: (a) structured classes with coordinate-aligned boundaries (e.g., hyper-rectangles and intersection of coordinate slabs), and (b) rotation-invariant or arbitrary convex geometries (e.g., Euclidean balls and general convex class). More precisely, typical collections of Borel subsets of $\mathbb{R}^d$ considered in literature include:
\begin{enumerate}[label=(\roman*)]
\item
$\mathcal{A}^{\mathrm{dist}} = \big\{ (-\infty, a_1] \times \dots \times (-\infty, a_d] : a_1, \dots, a_d \in \mathbb{R} \big\}$: the collection of all left-infinite rectangles (lower orthants),
\item 
$\mathcal{A}^{\max} = \big\{ (-\infty, t]^d : t \in \mathbb{R} \big\} = \big\{ \{y \in \mathbb{R}^d : \max_{1 \le j \le d} y_j \le t\} : t \in \mathbb{R} \big\}$: the collection of all left-infinite hypercubes,
\item
$\mathcal{A}^{\mathrm{rect}} = \big\{ \prod_{j=1}^d [a_j, b_j] : -\infty \le a_j \le b_j \le \infty \text{ for } j=1, \dots, d \big\}$: the collection of all hyperrectangles,
\item 
$\mathcal{A}^{\mathrm{hs}} = \big\{ \{y \in \mathbb{R}^d : \langle y, x \rangle \le t\} : x \in \mathbb{R}^d \setminus \{0\}, \, t \in \mathbb{R} \big\}$: the collection of all half-spaces,
\item 
$\mathcal{A}^{\mathrm{ball}} = \big\{ \{y \in \mathbb{R}^d : \|y - x\| \le r\} : x \in \mathbb{R}^d, \, r > 0 \big\}$: the collection of all Euclidean balls,
\item 
$\mathcal{A}^{\mathrm{con}}$: the collection of all Borel convex sets.
\end{enumerate}
Here and subsequently, $\|\cdot\|$ and $\langle \cdot, \cdot \rangle$ denote the Euclidean norm and standard inner product on $\mathbb{R}^d$, respectively. Clearly, $\mathcal{A}^{\mathrm{hs}}, \mathcal{A}^{\mathrm{ball}} \subset \mathcal{A}^{\mathrm{con}}$ and $\mathcal{A}^{\max} \subset \mathcal{A}^{\mathrm{dist}} \subset \mathcal{A}^{\mathrm{rect}} \subset \mathcal{A}^{\mathrm{con}}$. For any of the measurable classes $\mathcal{A}$ listed in (i)-(vi), the high-dimensional CLT is valid if the Kolmogorov–Smirnov-type distance satisfies
\begin{equation}
\rho_{n, \mathcal{A}} \equiv \sup_{A \in \mathcal{A}} \big| \mathbb{P}(T_n \in A) - \mathbb{P}(Z \in A) \big| \longrightarrow 0, \quad \text{as } n \to \infty,
\label{the CLT-1}
\end{equation}
under appropriate moment and regularity conditions. 
Early investigations into high dimensional CLT focused on general convex sets $\mathcal{A}^{\mathrm{con}}$ under finite third moments. \citet{senatov1981uniform} first showed that \eqref{the CLT-1} holds over compact convex sets when $d = o(n^{1/5})$. In a landmark paper, \citet{bentkus2003dependence} established that for independent and identically distributed summands with finite third moments, the Berry--Esseen bound on $\mathcal{A}^{\mathrm{con}}$ satisfies
\begin{equation}
    \rho_{n, \mathcal{A}^{\mathrm{con}}} \le C d^{1/4} \frac{\mathbb{E}\|X_1\|^3}{\sqrt{n}},
    \label{eq:bentkus-bound}
\end{equation}
which implies that $d = o(n^{2/7})$ is sufficient for CLT to hold. This bound was later extended to non-identically distributed summands with explicit constants by \citet{raic2019multivariate}. The dimension factor $d^{1/4}$ in \eqref{eq:bentkus-bound} stems from high-dimensional Gaussian 
surface area of a convex set in $\mathbb{R}^d$ (cf. \citet{nazarov2003maximal}). 
Recently, \citet{fang2024large} addressed this gap by establishing fourth-moment error bounds via Stein's method and exchangeable pairs. As they established, $\rho_{n,\mathcal{A}^{con}}\to 0$ as $n\to \infty$ if $d=o(n^{2/5})$ under finite fourth moments. For $\A = \A^{ball}$, the best result with respect to dimension dependence is due to \citet{fang2024large} who established (\ref{the CLT-1}) with $d = o(n^{1/2})$ under the existence of fourth moment of $X_1$.
In parallel, motivated by modern statistical applications in simultaneous inference and multiple testing, an active line of research has focused on coordinate-aligned classes such as hyper-rectangles $\mathcal{A}^{\mathrm{rect}}$ and $\mathcal{A}^{\max}$. In a series of foundational works, \citet{chernozhukov2013gaussian,chernozhukov2017central} established that the CLT remains valid over $\mathcal{A}^{\max}$ and $\mathcal{A}^{\mathrm{rect}}$ in ultra-high dimensions, permitting the dimension $d$ to grow as $\log d = o(n^{1/7})$ under sub-exponential moments, and conjectured that the optimal rate could be as fast as $\log d = o(n^{1/3})$. Subsequent investigations progressively relaxed these requirements (cf. \citet{lopes2020bootstrapping,koike2021notes,kuchibhotla2021high,chernozhukov2022improved}), using Gaussian anti-concentration on individual coordinates (cf. \citet{nazarov2003maximal}) to keep the error dependence on $d$ logarithmic rather than polynomial.

In this paper we focus on $\A^{hs}, \A^{ball}$ and $\A^{con}$ and we explore the optimal growth rate of the dimension $d$ as a function of the sample size $n$ for
the validity of the CLT in high dimensions under the existence of fourth moments. We observe that $d= o(\sqrt{n})$ is necessary for $\rho_{n, \A^{hs}}$ to converge to $0$ as $n\rightarrow \infty$. We utilized a construction of a sequence of iid mean zero random vectors with covariance matrix $I_d$ and $\sup_{n, d}\max_{i,j}\mathbb{E}|X_{ij}|^4 < \infty$. We translate the similar necessity for CLT over $\A^{ball}$ by essentially reducing $\rho_{n, \A^{ball}}$ to $\rho_{n, \A^{hs}}$ which is motivated by Vlasov's theorem (cf. \citet{vlasov1973approximative,BeauzamyMaurey1977}). On the other hand, we constructed the iid sequence with mean $0$, variance $I_d$ and $\sup_{n, d}\max_{i,j}d^{-3}\mathbb{E}|X_{ij}|^4 < \infty$ to claim that $d=o(n^{2/5})$ is necessary for $\rho_{n, \A^{con}}$ to converge to $0$ as $n\rightarrow \infty$. Moreover, if we relax the fourth moment condition to $\sup_{n, d}\max_{i,j}d^{-(3+l/2)}\mathbb{E}|X_{ij}|^4 < \infty$ for some $l \in \mathbb{N} \cup \{0\}$, then $d= o(n^{2/(5+l)})$ becomes necessary for the convex sets.
The rest of the paper is organized as follows. In sections \ref{sec:hs}, \ref{sec:ball} and \ref{sec:con}, we  provide the details of the necessity results for  $\mathcal{A}^{hs}$, $\mathcal{A}^{ball}$ and $\mathcal{A}^{con}$ respectively. In section \ref{sec:dis}, we discuss the necessary results of previous sections.

\section{Half Spaces}\label{sec:hs}
In this section we investigate the optimal growth of the dimension $d$ in the high dimensional CLT over half spaces. More specifically, we explore the growth of $d$ with respect to $n$ necessary to achieve
\begin{align}
    \rho_{n, \A^{hs}}\equiv   \sup_{A\in \A^{hs}} \Big|\mathbbm{P}(T_n\in A) -\mathbbm{P}(Z\in A)   \Big| \raw 0,  ~\mbox{as}~ \nti.
\end{align}
Since for any $a\in \mathbb{R}^d$, $\langle T_n, a\rangle$ is a scaled mean of univariate random variables (although depends on the choice of $a$), the bound on $\rho_{n, \A^{hs}}$ can be obtained by the classical Berry-Esseen theorem (cf. Theorem 12.4 of \citet{bhattacharya1986normal}) which entails 
\begin{align*}
    \rho_{n, \A^{hs}} \preceq\sup_{a \in \mathbb{R}^d}\dfrac{\sum_{i=1}^{n}\mathbb E|\langle a, X_i\rangle|^3}{n^{3/2}|a^\top \Sigma a|^{3/2}}&\leq \sup_{a \in \mathbb{R}^d}\dfrac{\|a\|^2\sum_{i=1}^{n}\Big\{\mathbb E|\langle a, X_i\rangle|^2\mathbb E\|X_i\|^4\Big\}}{n^{3/2}|a^\top \Sigma a|^{3/2}}^{1/2}\\
    & \leq \sup_{a \in \mathbb{R}^d}\dfrac{\|a\|^2\Big(\sum_{i=1}^{n}\mathbb E\|X_i\|^4\Big)}{n|a^\top \Sigma a|}^{1/2}\\
    & \leq\|\Sigma^{-1}\|_{op}\Big(n^{-2}\sum_{i=1}^{n}\mathbb E\|X_i\|^4\Big)^{1/2},
\end{align*}
To obtain either of second, third and fourth inequality, we have only utilized the Cauchy-Schwarz inequality.  Note that the Berry-Esseen rate in this class matches with the Berry-Esseen rate for $\mathcal{A}^{ball}$, obtained recently by \citet{fang2024large} up to a logartithmic factor. If the eigen values of $\Sigma$ are bounded away from $0$ and $\sup_{n, d}\max_{1\leq i \leq n, 1\leq j \leq d}\mathbb E|X_{ij}|^4 \preceq 1$ then $$n^{-2}\sum_{i=1}^{n}\mathbb E\|X_i\|^4 \preceq d^2/n.$$ On the other hand, since there exists $b\in \mathbb{R}^d$ such that $a = \Sigma^{-1/2}b$ due to invertibility of $\Sigma$, one can replace $a$ by $\sigma^{-1/2}b$ in the calculations of Berry-Esseen rate obtained above. In that way one can also have $$\rho_{n, \A^{hs}} \preceq \Big(n^{-2}\sum_{i=1}^{n}\mathbb E\|\Sigma^{-1/2}X_i\|^4\Big)^{1/2}.$$  If $X_1,\dots, X_n$ are mean zero independent and identically distributed (iid) random vectors then $$n^{-2}\sum_{i=1}^{n}\mathbb E\|\Sigma^{-1/2}X_i\|^4 \geq n^{-2}\sum_{i=1}^{n}(\mathbb E\|\Sigma^{-1/2}X_i\|^2)^2= d^2/n.$$ Thus the typical order of the high dimensional Berry-Esseen rate is $d/\sqrt{n}$, under the existence of fourth moments. 
Therefore, in general, the high dimensioonal CLT over half spaces hold when $d = o(\sqrt{n})$. In this paper, we would like to ask whether $d=o(\sqrt{n})$ is also necessary for $\rho_{n, \A^{hs}} \rightarrow 0$ under existence of fourth moments. We have the following theorem.
\begin{theorem}\label{thm:hs}
    Let $X_1,\dots, X_n$ be iid random vectors. Then there exists a $d$-dimensional distribution of $X_1$ with $\mathbb EX_1 = 0$, $Var(X_1) = I_d$ and $\sup_{n, d}\max_{ 1\leq j \leq d}\mathbb E|X_{1j}|^4 \preceq 1$ such that if $ \rho_{n, \A^{hs}} \rightarrow 0$ as $d, n \rightarrow \infty$ then we must have $d = o(\sqrt{n})$. 
\end{theorem}
Proof: Note that if $\rho_{n, \A^{hs}} \rightarrow 0$ as $d, n \rightarrow \infty$, then by continuity of a probability measure, we must have
$\mathbbm{P}(\langle T_n, a \rangle = t)) \rightarrow 0$, as $d, n \rightarrow \infty$, for any $a$ and $t$. Thus it is enough to show that
\begin{align}\label{eqn:2.1.1}
    |\mathbbm{P}(\langle T_n, u \rangle = t_0)| \succeq e^{-n/d^2},
\end{align}
for some particular choice of $t_0$ and a unit vector $u$ which we specify later. To conclude (\ref{eqn:2.1.1}), we have to construct the iid random vectors $X_1,\dots, X_n$ (the distribution of which may depend on $n$ and $d$) to be such that
\begin{enumerate}[label = (\alph*)]
\item the one dimensional distribution of $\langle n^{-1/2}X_1, u \rangle$ has a jump of order $e^{-1/d^2}$. 
\item $\mathbb EX_1 = 0$ and $Var(X_1) = I_d$.
    \item $\sup_{n, d}\max_{ 1\leq j \leq d}\mathbb E|X_{1j}|^4 \preceq 1$.
\end{enumerate}
Let $B_{1}, \dots, B_{n}$ be iid with $B_{1}\sim \text{Bernoulli}(\pi_d)$ where $\pi_d = (1+d^2)^{-1}$. For each $i \in \{1,\dots, n\}$, set $C_{ i} = [\pi_d(1-\pi_d)]^{-1/2}(B_{i} - \pi_d)$, the centered and scaled version of $B_i$. Suppose that $U_1,\dots U_n$ are iid random vectors from $N(0, I_d)$ and assume them to be independent with $B_{1},\dots, B_n$. Finally for any $i\in \{1,\dots, n\}$, define
$$X_i = C_i u + (I_d-uu^\top) U_i\;\;\text{and}\;\; T_n:=n^{-1/2}\sum_{i=1}^n X_i.$$
Note that, $\langle T_n,u\rangle=n^{-1/2}\sum_{i=1}^n C_i=[n\pi_d(1-\pi_d)]^{-1/2}(\sum_{i=1}^nB_{i} - n\pi_d)$ with $\sum_{i=1}^nB_{i}\sim \text{Binomial}(n,\pi_d)$ and  $$\langle n^{-1/2}X_i, u \rangle = n^{-1/2}C_i= \dfrac{B_{i} - \pi_d}{\sqrt{n\pi_d(1-\pi_d)}}.$$ Now without loss of generality choose, $t_0 = -n\pi_d/\sqrt{n\pi_d(1-\pi_d)}$. Therefore,
\begin{align*}
\mathbb{P}[\langle T_n,u\rangle=t_0]=\mathbb{P}\Big[\sum_{i=1}^nB_i=0\Big]=\Big[e^{\log(1-(1+d^2)^{-1})}\Big]^n\geq e^{-n/d^2},    \end{align*}
where we use $\log (1-x)\geq -x(1-x)^{-1}$ for $x \in (0, 1)$. Thus we have (a). Again, (b) is trivially true by the construction. It remains to validate (c). 
For that, let $X_{1j}$ denote the $j$-th coordinate of $X_1$ and $e_j$ denote the $j$-th canonical basis for $j \in \{1, \dots, d\}$. Therefore we can have,
\begin{align}\label{eqn:2.1.2}
X_{1j}=\langle e_j, X_1\rangle=C_1u_j+\langle e_j, (I_d - uu^\top) U_1\rangle \stackrel{d}{=} C_1 u_j + \sqrt{1 - u_j^2} \, N_j,
\end{align}
where $N_1, \dots, N_d$ are distributed as $N(0, 1)$ and independent of $C_1$ by construction.
Since $B_{i} \sim \text{Bernoulli}(\pi_d)$, we have
\[
\mathbb{E}[(B_{i} - \pi_d)^4] = \pi_d(1-\pi_d) \left(1 - 3\pi_d(1-\pi_d)\right) \le \pi_d(1-\pi_d).
\]
Recalling $\pi_d = \frac{1}{d^2+1}$, we have $\pi_d(1-\pi_d) = \frac{d^2}{(d^2+1)^2}$ and hence
\[
\mathbb{E}[C_{i}^4] = \frac{\mathbb{E}[(B_{i} - \pi_d)^4]}{(\pi_d(1-\pi_d))^2} \le \frac{1}{\pi_d(1-\pi_d)} \le \frac{(d^2+1)^2}{d^2} \le 4d^2.
\]
Again we have
$\mathbb{E}[N_j^4] = 3 $. Thus from (\ref{eqn:2.1.2}), for any $j\in \{1,\dots, d\}$ we have
\[\mathbb{E}|X_{1j}|^4 \leq 8\Big(4u_j^4d^2 + 3(1-u_j^2)^2\Big) \preceq 1,\]
provided $u_j = 1/\sqrt{d}$. Thus we have (c) if we set $u = d^{-1/2}(1,\dots, 1)^\top$.
Therefore the proof is complete.\hfill $\square$

\section{Euclidean Balls}\label{sec:ball}
In this section, we explore the optimal growth rate of the dimension $d$ in the high dimensional CLT over Euclidean balls. We investigate the optimal growth of $d$ (with respect to $n$) above which
\begin{align}
    \rho_{n, \A^{ball}}\equiv   \sup_{A\in \A^{ball}} \Big|\mathbbm{P}(T_n\in A) -\mathbbm{P}(Z\in A)   \Big| \raw 0,  ~\mbox{as}~ \nti.
\end{align}
fails. \citet{fang2024large} recently established under existence of fourth moments that
$$\rho_{n, \mathcal{A}^{ball}}\preceq \|\Sigma^{-1}\|_{op}\Big(n^{-2}\sum_{i=1}^{n}\mathbb E\|X_i\|^4\Big)^{1/2}. $$ This high dimensional Berry-Esseen rate matches with the Berry-Esseen rate for the half spaces, derived in the previous section. As mentioned earlier, the typical order of the above Berry-Esseen rate is $d/\sqrt{n}$. Thus $\rho_{n, \A^{ball}} \rightarrow 0$, as $n\rightarrow \infty$ when $d = o(\sqrt{n})$. In the next theorem we establish that $d = o(\sqrt{n})$ is also necessary for $\rho_{n, \A^{ball}}\rightarrow 0$, under the existence of fourth moments. We would like to point out that \citet{fang2024large} established that $d=o(n)$ is both necessary and sufficient for high dimensional CLT uniformly over centered Euclidean balls. 
\begin{theorem}\label{thm:ball}
    Let $X_1,\dots, X_n$ be iid random vectors. Then there exists a $d$-dimensional distribution of $X_1$ with $\mathbb EX_1 = 0$, $Var(X_1) = I_d$ and $\sup_{n, d}\max_{ 1\leq j \leq d}\mathbb E|X_{1j}|^4 \preceq 1$ such that if $ \rho_{n, \A^{ball}} \rightarrow 0$ as $d, n \rightarrow \infty$ then we must have $d = o(\sqrt{n})$. 
\end{theorem}

Proof: It is enough to show that $\rho_{n, \A^{hs}} \leq \rho_{n, \A^{ball}}$, since then Theorem \ref{thm:hs} will give us the result. This idea is motivated by the Vlasov's theorem (cf. \citet{vlasov1973approximative,BeauzamyMaurey1977}). 
First fix a $d$-dimensional unit vector $u$ and a real number $t$. Then consider a sequence of Euclidean balls $\{B_m(u, t)\}_{m\geq 1}$ as
\[
B_m(u, t) := \left\{ \bm{y} \in \mathbb{R}^d : \|\bm{y} + mu\|^2 \le m^2 + 2tm-1 \right\}.
\]
If $m^2 + 2tm-1 < 0$, then the corresponding set is identified as the null set. Now note that for any $\bm{z}\in\mathbb R^d$,
\[
\bm{z} \in B_m(u, t) \iff u^\top \bm{z} \le t - \frac{\|\bm{z}\|^2+1}{2m}.
\]
Clearly, $\{B_m(u,t)\}_{m \geq 1}$ is a non-decreasing sequence of sets and hence  
\(\lim_{m \uparrow \infty} \mathbb{P}(V \in B_m(u, t)) = \mathbb{P}\left(V \in \bigcup_{m \geq 1} B_m(u, t) \right),\) for any $d$-dimensional random vector $V$.
 Now, we claim that:
\[
\bigcup_{m\geq 1} B_m(t) = \left\{y\in\mathbb R^d: \langle y,u\rangle< t \right\}.
\]
For the $(\supseteq)$ side, fix $z\in \left\{y\in\mathbb R^d: \langle y,u\rangle< t \right\}$. Then we have $\langle z,u\rangle\le t-\varepsilon,\;\text{for some}\;\varepsilon>0.$ 
    Since 
    $\frac{\|z\|^2+1}{2m}\le \varepsilon$ for large enough $m$, $z \in B_m(u, t)$ for large enough $m$.
 Conversely for $(\subseteq)$ side, fix $z \in \bigcup_{m \geq 1} B_m(u, t)$. Then there exists some natural number $m_0$ such that
    $\langle z,u\rangle \le t - \frac{\|\bm{z}\|^2+1}{2m_0}$ and hence $\langle z,u\rangle < t$.

Based on all the above we have for any unit vector $u \in \mathbb{R}^d$ and $t\in \mathbb{R}$,
\begin{align*}
\rho_{n, \A^{ball}} &= \sup_{B \in \mathcal{A}^{ball}} |\mathbb{P}(T_{n} \in B) - \mathbb{P}(Z\in B)|\\
&\ge \lim_{m\rightarrow \infty}|\mathbb{P}(T_{n} \in B_m(u, t)) - \mathbb{P}\left(Z \in B_m(u, t)\right)|\\
&= \Big|\mathbb{P}\Big(T_{n} \in \bigcup_{m=1}^{\infty}B_m(u, t)\Big) - \mathbb{P}\Big(Z \in \bigcup_{m=1}^{\infty} B_m(u, t)\Big)\Big|\\
&= \Big|\mathbb{P}\big(\langle T_{n}, u\rangle < t\big) - \mathbb{P}\big(\langle Z, u\rangle < t\big)\Big|
\end{align*}
Thus by right continuity of the distribution function of $\langle T_{n}, u\rangle$, we have
\begin{align*}
\rho_{n, \A^{ball}} & \geq \sup_{t \in \mathbb{R},\|u\|=1,} \Big|\mathbb{P}\big(\langle T_{n}, u\rangle \leq t\big) - \mathbb{P}\big(\langle Z, u\rangle \leq t\big)\Big| = \rho_{n, \A^{hs}}.
\end{align*}
Therefore, the proof is complete. \hfill $\square$

\section{Convex sets}\label{sec:con}

In this section, we establish the optimal growth rate of the dimension $d$ in the high dimensional CLT over the class of Borel convex sets. Equivalently, we seek for the  optimal rate of growth of $d$ (with respect to $n$) above which
\begin{align}
    \rho_{n, \A^{con}}\equiv   \sup_{A\in \A^{con}} \Big|\mathbb{P}(T_n\in A) -\mathbb{P}(Z\in A)   \Big| \raw 0,  ~\mbox{as}~ \nti.
\end{align}
fails. \citet{fang2024large} established under existence of fourth moments that
$$\rho_{n, \mathcal{A}^{ball}}\preceq d^{1/4}\Big(n^{-2}\sum_{i=1}^{n}\mathbb E\|\Sigma^{-1/2}X_i\|^4\Big)^{1/2}, $$ up to a logarithmic factor. The distinction of the factor $d^{1/4}$ in case of class of convex sets and Euclidean Balls, arises due to the respective Gaussian isoperimetric bounds for these two classes (cf. \citet{bentkus2003dependence,nazarov2003maximal,zhilova2020nonclassical}). Since the term inside the braces in the RHS of the above display is in general of order $d^2/n$, $\rho_{n, \A^{con}} \rightarrow 0$, as $n\rightarrow \infty$ when $d = o(n^{2/5})$. Here we explore whether $d = o(n^{2/5})$ is also necessary under existence of fourth moments. 
We have the following theorem.
\begin{theorem}\label{thm:con}
    Let $X_1,\dots, X_n$ be iid random vectors. Then there exists a $d$-dimensional distribution of $X_1$ with $\mathbb EX_1 = 0$, $Var(X_1) = I_d$ and $\sup_{n, d}\max_{ 1\leq j \leq d}d^{-3}\mathbb E|X_{1j}|^4 \preceq 1$ such that if $ \rho_{n, \A^{con}} \rightarrow 0$ as $d, n \rightarrow \infty$ then we must have $d = o(n^{2/5})$. 
\end{theorem}

Proof: First note that for any two Borel convex sets $C_1\subseteq C_2$, 
\[   2\sup_{C \in \mathcal{A}^{con}} \big|\mathbb{P}(T_{n} \in C) - \mathbb{P}(Z \in C)\big| 
    \ge  \mathbb{P}\big(Z \in C_2 \setminus C_1\big) - \mathbb{P}\big(T_n \in C_2 \setminus C_1\big) .\]
We are going to show that for some choices of $C_1$ and $C_2$, the first term in the RHS of above display is $\succeq \big(1-e^{-\big(d^{5/4}/n^{1/2}\big)}\big)$ while a suitable choice of the distribution of $X_1$ gives the second term to be exactly $0$. Thus we must have $d = o(n^{2/5})$ for LHS to converge to $0$ as $n, d \rightarrow \infty$. To get the aforementioned lower bound, first let $\rho = d^{1/4}$, $\delta = d/\sqrt{n}$ and 
\begin{equation}\label{eq:def_N}
    N := \left\lfloor c_0 \rho \exp\left\{\frac{\rho^2}{2}\right\} \right\rfloor,
\end{equation}
where $c_0 = 9\sqrt{2\pi}/(8 e^3)$. 
Subsequently define $\mathcal{Z} := \big(\{-1/\sqrt{d}, \, 1/\sqrt{d}\}^d\big)^N$. 
If $\bm{z} = (z_1, \dots, z_N) \in \mathcal{Z}$, then each $z_k$ satisfies $\|z_k\|^2 = \sum_{j=1}^d (\pm 1/\sqrt{d})^2 = 1$ and hence each $z_k$ lies on $\mathbb{S}^{d-1}$. For any such $\bm{z} \in \mathcal{Z}$, define the nested convex polyhedra $A_1(\bm{z}) \subset A_2(\bm{z})$ by:
\begin{align}\label{eq:polyhedra_def}
    A_1(\bm{z}) = A(\rho; \bm{z}) &:= \bigcap_{k=1}^N \big\{x \in \mathbb{R}^d : \langle x, z_k \rangle \le \rho\big\}, \nonumber \\
    A_2(\bm{z}) = A(\rho + \delta; \bm{z}) &:= \bigcap_{k=1}^N \big\{x \in \mathbb{R}^d : \langle x, z_k \rangle \le \rho + \delta\big\}.
\end{align}
Note that, we can assume that $\delta = o(1)$ since $d = o(\sqrt{n})$ is necessary for $\rho_{n,\A^{hs}} \rightarrow 0$, as $n, d \rightarrow \infty$, and $\rho_{n,\A^{hs}} \leq \rho_{n,\A^{con}}$.

We first show that 
 \begin{align}\label{eqn:4.100}
\mathbb{P}\big(Z \in A_2(z^*) \setminus A_1(z^*)\big) \succeq \big(1-e^{-\big(d^{5/4}/n^{1/2}\big)}\big),
\end{align}
for some $z^* \in \mathcal{Z}$. Hence from now on we can consider $C_i = A_i(z^*)$, $i=1, 2$. We will
utilize the random set construction in Lemma 2.6 of \citet{bentkus2003dependence}). To that end, let $Z_1, \dots, Z_N$ be iid random vectors, each distributed uniformly on $\{-1/\sqrt{d}, 1/\sqrt{d}\}^d$ and  defined by:
\[
    \sqrt{d} Z_k = (\varepsilon_{k1}, \dots, \varepsilon_{kd})^\top, \qquad k \in \{1, \dots, N\},
\]
where $\{\varepsilon_{kj} : 1 \le k \le N, 1 \le j \le d\}$ are i.i.d.\ Rademacher variables. We denote the random polyhedra by $A_1 := A_1(\bm{Z})$ and $A_2 := A_2(\bm{Z})$, and the random shell volume by $W := \mathbb{P}(Z \in (A_2 \setminus A_1))$ and the law of $N(0, I_d)$ by $\gamma_d$. Thus to claim (\ref{eqn:4.100}), it is enough to show that 
\begin{align}\label{eqn:4.103}
\mathbb{E}_{\bm{Z}}(W)\succeq \big(1-e^{-\big(d^{5/4}/n^{1/2}\big)}\big).
\end{align}
Now by equation (2.29) of \citet{bentkus2003dependence}), 
it follows that
\begin{equation}\label{eq:J_integral}
\mathbb{E}_{\bm{Z}}[W]\ge e^{-1/2}NJ,
\qquad
J:=\int_{\substack{\|x\|\le\sqrt d\\
\rho<x_1\le\rho+\delta}}\gamma_d(dx).
\end{equation}

Now by Fubini's theorem we have
\begin{align}\label{eqn:4.101}
    J = \int_\rho^{\rho + \delta} \frac{1}{\sqrt{2\pi}} e^{-x_1^2/2} \left[ \int_{\sum_{j=2}^d x_j^2 \le d - x_1^2} \gamma_{d-1}(dx') \right] dx_1.
\end{align}
Since $x_1 \in [\rho, \rho + \delta]$ and $\rho = d^{1/4}$, we have $x_1 \le d^{1/4} + \delta$ and hence 
\[
    d - x_1^2 \ge d - \sqrt{d} - O\big(\delta d^{1/4}\big).
\]
Thus with some $\eta_2, \dots, \eta_d \stackrel{\text{i.i.d.}}{\sim} \text{N}(0, 1)$, we have
\begin{align*}
    \int_{\sum_{j=2}^d x_j^2 \le d - x_1^2} \gamma_{d-1}(dx') & = \mathbb{P}\left(\sum_{j=2}^d \eta_j^2 \le d - x_1^2\right) \\
    & \ge \mathbb{P}\left(\frac{\sum_{j=2}^d (\eta_j^2 - 1)}{\sqrt{2(d-1)}} \le \frac{1 - \sqrt{d} - O(\delta d^{1/4})}{\sqrt{2(d-1)}}\right).
\end{align*}
Now note that
\[
    \frac{-\sqrt{d} + 1 - O(\delta d^{1/4})}{\sqrt{2}\sqrt{d}} 
    = -\frac{1}{\sqrt{2}} + \frac{1}{\sqrt{2d}} - O\left(\frac{\delta}{d^{1/4}}\right) \longrightarrow -\frac{1}{\sqrt{2}} \approx -0.7071 \quad \text{as } d \to \infty,
\]
since $\delta / d^{1/4} = o(1)$ due to $\delta = o(1)$. Thus by the Central Limit Theorem we have 
\[
    \int_{\sum_{j=2}^d x_j^2 \le d - x_1^2} \gamma_{d-1}(dx') \ge \frac{1}{2} \Phi(-1/\sqrt{2}) > 0, \qquad \text{for all } x_1 \in [\rho, \rho + \delta].
\]
Thus from (\ref{eqn:4.101}) we have
\begin{align}\label{eqn:4.102}
    J \ge \Phi(-1/\sqrt{2})(8 \pi)^{-1/2} \int_\rho^{\rho + \delta} e^{-x_1^2/2} dx_1.
\end{align}
Applying the linear change of variables $x_1 = \rho + u$ with $u \in [0, \delta]$, we have 
\[
    \int_\rho^{\rho + \delta} e^{-x_1^2/2} dx_1 
    = e^{-\rho^2/2} \int_0^\delta e^{-\rho u - u^2/2} du.
\]
Because $u \in [0, \delta]$ and $\delta = o(1)$, we have $u^2/2 \le \delta^2/2 \to 0$; in particular, $e^{-u^2/2} \ge e^{-1/2} \ge 1/2$ for all sufficiently large $n$. Thus:
\[
    \int_0^\delta e^{-\rho u - u^2/2} du 
    \ge \frac{1}{2} \int_0^\delta e^{-\rho u} du 
    = \frac{1}{2} \left[ \frac{-e^{-\rho u}}{\rho} \right]_0^\delta 
    = \frac{1}{2} \left( \frac{1 - e^{-\rho \delta}}{\rho} \right).
\]
Combining above steps, (\ref{eqn:4.102}) reduces to
\begin{equation}\label{eq:J_lower}
    J \ge \Phi(-1/\sqrt{2})(8 \pi)^{-1/2} e^{-\rho^2/2} \left(\frac{1 - e^{-\rho \delta}}{\rho}\right).
\end{equation}
Substituting \eqref{eq:J_lower} into \eqref{eq:J_integral} and using $N \ge \frac{c_0}{2} d^{1/4} e^{\rho^2/2}$ from \eqref{eq:def_N}:
\begin{align}\label{eq:EW_evaluated}
    \mathbb{E}_{\bm{Z}}[W] 
    &\ge e^{-1/2} \left(\frac{c_0}{2} d^{1/4} e^{\rho^2/2}\right)  \Phi(-1/\sqrt{2})(8 \pi)^{-1/2} e^{-\rho^2/2} \left(\frac{1 - e^{-\rho \delta}}{\rho}\right) \nonumber \\
    &= \alpha  \big(1 - e^{-\rho \delta}\big),
\end{align}
where $\alpha = c_0e^{-1/2} \Phi(-1/\sqrt{2})(32 \pi)^{-1/2}$ with $c_0 = 9\sqrt{2\pi}/(8e^3)$. Thus we establish (\ref{eqn:4.103}), and hence
\begin{align}\label{eqn:4.105}
    \mathbb{P}\big(Z \in C_2 \setminus C_1\big) \succeq (1-e^{-(d^{5/4}/n^{1/4})}).
\end{align}

Now we move towards construction of the distribution of $X_1$ (and hence of $T_n$) so that we obtain $\mathbb{P}(T_n \in C_2 \setminus C_1) = 0$.  
Let $X_{1}, \dots, X_{n}$ be iid random vectors in $\mathbb{R}^d$ with $X_1$ have iid coordinates. Let $X_{11}$ has the symmetric three-point distribution:
\begin{equation}\label{eq:Xij_dist}
    X_{11} = \begin{cases} 
        0, & \text{with probability } 1 - 2p_d, \\
        -c, & \text{with probability } p_d, \\
        +c, & \text{with probability } p_d,
    \end{cases}
    \qquad \text{where } p_d := \frac{1}{2c^2},
\end{equation}
for a jump size $c > 0$ specified below. Note that:
\[
    \mathbb{E}[X_{11}] = 0, \qquad \mathbb{E}[X_{11}^2] = 2p_d c^2 = 1, \qquad \mathbb{E}[X_{11}^4] = 2p_d c^4 = c^2,
\]
implying that $\mathbb{E}[X_{1}] = \mathbf{0}$ and $\operatorname{Cov}(X_{1}) = I_d$. The normalized sum is defined by $T_{n} := \frac{1}{\sqrt{n}} \sum_{i=1}^n X_{i}$.
Recall that each deterministic $z_k^* \in \mathcal{Z} = \big(\{-1/\sqrt{d}, 1/\sqrt{d}\}^d\big)^N$ has the form:
\[
    z_k^* = \frac{1}{\sqrt{d}} \left(\varepsilon_{k1}^*, \dots, \varepsilon_{kd}^*\right)^\top, \qquad \varepsilon_{kj}^* \in \{-1, +1\}.
\]
Projecting the normalized sum $T_{n}$ onto $z_k^*$ yields:
\begin{align}\label{eq:projection_expansion}
    \langle T_{n}, z_k^* \rangle 
    &= \frac{1}{\sqrt{n}} \sum_{i=1}^n \langle X_{i}, z_k^* \rangle= \frac{1}{\sqrt{nd}} \sum_{i=1}^n \sum_{j=1}^d \varepsilon_{kj}^* X_{ij}= \frac{c}{\sqrt{nd}} \sum_{i=1}^n \sum_{j=1}^d \varepsilon_{kj}^* \left(\frac{X_{ij}}{c}\right).
\end{align}
Because $X_{ij} \in \{-c, 0, c\}$, the normalized variable $X_{ij}/c$ takes values in $\{-1, 0, +1\}$. Since $\varepsilon_{kj}^* \in \{-1, +1\}$, the product $\varepsilon_{kj}^*(X_{ij}/c)$ is an integer for all $i, j$. Consequently, 
for every $k \in \{1, \dots, N\}$, $$\langle T_{n}, z_k^* \rangle \in h\mathbb{Z} \; \text{with probability 1, where}\; h := \frac{c}{\sqrt{nd}}.$$ To ensure $\mathbb{P}(T_n \in C_2 \setminus C_1) = 0$, we need to make sure that no points of $h\mathbb{Z}$ fall inside the half-open interval $(\rho, \rho + \delta]$. 
For that we require to have an integer $m$ such that
\begin{equation}\label{eq:grid_conditions}
    \rho = m h  \qquad \text{and} \qquad h > \delta.
\end{equation}
Because $h > \delta$, 
\begin{equation}\label{eq:empty_intersection}
    (\rho, \rho + \delta] \subset \big(mh, \, (m+1)h\big) \implies (\rho, \rho + \delta] \cap h\mathbb{Z} = \varnothing.
\end{equation}
Now to satisfy \eqref{eq:grid_conditions}, we require $m = \frac{\rho}{h} < \frac{\rho}{\delta}$. Since $\frac{\rho}{\delta} = \frac{d^{1/4}}{d n^{-1/2}} = \frac{\sqrt{n}}{d^{3/4}} \to \infty$ as $n \to \infty$ (since $d=o(\sqrt{n})$), we set:
\begin{equation}\label{eq:m_choice}
    m := \left\lfloor \frac{\rho}{2\delta} \right\rfloor = \left\lfloor \frac{\sqrt{n}}{2 d^{3/4}} \right\rfloor \in \mathbb{N},\;\;\;\;\; h := \frac{\rho}{m} \;\;\;\;\; \text{and}\;\;\;\; c := h\sqrt{nd}.
\end{equation}
With this choice of $c$ we have $\mathbb{P}(T_n \in C_2 \setminus C_1) = 0$ due to (\ref{eq:empty_intersection}). Furthermore, for all sufficiently large $n$, we have $m \ge \frac{\rho}{2\delta} - 1 \ge \frac{\rho}{4\delta}$, which yields:
\[
    h \le 4\delta = \frac{4d}{\sqrt{n}} \implies c = h\sqrt{nd} \le \left(\frac{4d}{\sqrt{n}}\right) \sqrt{nd} = 4 d^{3/2}. 
\]
Consequently, there exist $d-$dimensional iid random vectors $X_1,\dots, X_n$ such that $\mathbb{E}[X_{1}] = \mathbf{0}$, $\operatorname{Cov}(X_{1}) = I_d$, 
$\sup_{n,d} \sup_{i,j}d^{-3} \mathbb{E}[X_{ij}^4] = d^{-3}c^2 \le 16.$ and 
\begin{align}\label{eqn:4.106}
    \mathbb{P}\big(n^{-1/2}\sum_{i=1}^{n}X_i \in C_2 \setminus C_1\big)= 0.
\end{align}
Combining (\ref{eqn:4.105}) and (\ref{eqn:4.106}), the proof is now complete. $\hfill \square$

\section{Discussions}\label{sec:dis}
In this paper, we explore the dimension growth rate necessary for high dimensional CLT to hold uniformly over half spaces or over Euclidean balls or over Borel convex sets.  In section \ref{sec:hs}, we show that $d= o(\sqrt{n})$ is necessary for high dimensional CLT over half spaces. In Section \ref{sec:ball}, we extend the arguments of \ref{sec:hs} to establish that $d= o(\sqrt{n})$ is also necessary for high dimensional CLT over Euclidean balls. In section \ref{sec:con}, we presumed $d=o(\sqrt{n})$ and utilized it to show that $d=o(n^{2/5})$ is necessary for CLT over Borel convex sets. In each of theorems \ref{thm:hs}, \ref{thm:ball} and \ref{thm:con}, we considered random vectors $X_1,\dots, X_n$ to be iid with mean $0$, variance $I_d$. However, there is clear distinction between the construction corresponding to Theorems \ref{thm:hs} and \ref{thm:ball} and that corresponding to Theorem \ref{thm:con}. The construction corresponding to Theorem~\ref{thm:hs} and \ref{thm:ball} permits coordinate random variables with uniformly bounded fourth moments (i.e. $\sup_{n, d} \max_{1\leq i \leq n, 1 \le j \le d} \mathbb{E}[X_{ij}^4] \preceq 1$), whereas the  construction corresponding to the general convex sets in Theorem~\ref{thm:con} requires the fourth moment to scale as $\max_{1\leq i \leq n, 1 \le j \le d} \mathbb{E}[X_{ij}^4] \asymp d^3$. 

While Theorem \ref{thm:con} demonstrates that, the rate $d=o(n^{2/5})$ is necessary for CLT over class of convex sets to hold, the underlying  construction of the distribution of $T_n$ reveals an intriguing insight. 
Note that we started with $\delta=\frac{d}{\sqrt{n}}$ since $d=o(\sqrt{n})$ was observed to be necessary for CLT over half-spaces or over Euclidean balls. This results in the necessity of $d^{5/4}=o(\sqrt{n})$. Thus in the next loop, we could start with $\delta=\frac{d^{5/4}}{\sqrt{n}}$ in the proof of Theorem \ref{thm:con} and then we will get
\[
 \sup_{C \in \mathcal{A}^{con}} \big|\mathbb{P}(T_{n} \in C) - \gamma_d(C)\big| \ge \frac{\alpha}{2} \left( 1 - \exp\left\{ -\frac{d^{6/4}}{\sqrt{n}} \right\} \right).
\]
Moreover, following the lattice random vector construction, we can also get $\mathbb{P}(T_n \in C_2\setminus C_1) = 0$, however with  
$\sup_{n,d} \sup_{i,j}d^{-7/2} \mathbb{E}[X_{ij}^4] \preceq 1$. Thus essentially, we can claim that $d^{6/4} = o(\sqrt{n})$ is also necessary for CLT to hold over Borel convex sets. If we proceed this way then 
we actually have the following proposition.
\begin{proposition}
 Let $X_1,\dots, X_n$ be iid random vectors and fix $\ell\in\mathbb{N}\cup\{0\}$. Then there exists a $d$-dimensional distribution of $X_1$ (which may depend on $n$ and $d$) with $\mathbb EX_1 = 0$, $Var(X_1) = I_d$ and $$\sup_{n, d}\max_{ 1\leq j \leq d}d^{-(3+\ell/2)}\mathbb E|X_{1j}|^4 \preceq 1$$ such that if $ \rho_{n, \A^{con}} \rightarrow 0$ as $d, n \rightarrow \infty$ then we must have $d^{\frac{5+\ell}{4}} = o(\sqrt{n})$.
\end{proposition}
This behavior is  primarily due to the fact that we allow the fourth moment to depend on $d$. If we consider the fourth moments to be uniformly bounded, like in case of half spaces and Euclidean balls, then we believe that $d= o(n^{2/5})$ is the optimal rate possible for CLT over Borel convex sets. We conjecture that: 
\begin{quote}
  ``Let $X_1,\dots, X_n$ be iid random vectors. Then there exists a $d$-dimensional distribution of $X_1$ with $\mathbb EX_1 = 0$, $Var(X_1) = I_d$ and $\sup_{n, d}\max_{ 1\leq j \leq d}\mathbb E|X_{1j}|^4 \preceq 1$ such that if $ \rho_{n, \A^{con}} \rightarrow 0$ as $d, n \rightarrow \infty$ then we must have $d = o(n^{2/5})$."  
\end{quote}

\bibliographystyle{amsplain}

\end{document}